
\magnification\magstep1
\documentstyle{amsppt}
\TagsOnRight
\NoBlackBoxes
\hoffset=.1in

\define\sh{\operatorname{sh}}
\define\ch{\operatorname{ch}}

\define\ld{\ldots}

\define\a{\alpha}
\define\be{\beta}
\define\vf{\varphi}
\define\lm{\lambda}
\define\ve{\varepsilon}
\define\gm{\gamma}
\define\de{\delta}
\define\si{\sigma}
\define\om{\omega}

\define\vd{\varDelta}
\define\vl{\varLambda}
\define\vg{\varGamma}
\define\vo{\varOmega}

\define\vt{\varTheta}

\define\cam{\Cal M}
\define\can{\Cal N}

\define\car{\Cal R}

\define\xti{\tilde x}
\define\yti{\tilde y}

\define\pd#1#2{\dfrac{\partial#1}{\partial#2}}

\define\vc#1{(#1_1,\ldots,#1_n)}
\define\vct#1{[#1_1,\ldots,#1_n]}
\define\vect#1{\{#1_1,\ldots,#1_n\}}


\define\jak{\dfrac{i}{\kappa}}
\define\mak{\Big(\dfrac{P_0}{\kappa}\Big)}
\define\maks{\frac{P_0}{\kappa}}
\define\kap{{\Cal P}_\kappa}
\define\kapti{\widetilde{\Cal P}_\kappa}

\define\kub{\dfrac{2\kappa}}
\define\rhoti{\widetilde \rho}


\define\1{$\kappa$-Poincar\'e group}
\define\2{bicovariant}
\define\3{calculus}
\define\4{differential}
\define\5{dimensional}
\define\6{quantum}
\define\7{coefficient}
\define\8{operator}
\define\9{properties}
\define\0{algebra}

\define\kapo{$\kappa$-Poincar\'e}
\define\mink{$\kappa$-Min\-ko\-wski}
\define\linv{left-invariant}
\define\rinv{right-invariant}

\topmatter
\title The bicovariant differential calculus\\
on the $\kappa$-Poincar\'e group\\
and on the $\kappa$-Minkowski space
\endtitle
\rightheadtext{Differential calculi}
\author Piotr \/ Kosi\'nski*\\
{\it Department of Theoretical Physics}\\
{\it University of \L \'od\'z}\\
{\it ul. Pomorska 149/153, 90--236 \L \'od\'z, Poland}\\
Pawe\l \/ Ma\'slanka*\\
{\it Department of Functional Analysis}\\
{\it University of \L \'od\'z}\\
{\it ul. St. Banacha 22, 90--238 \L \'od\'z, Poland}\\
Jan \/ Sobczyk**\\
{\it Institute of Theoretical Physics}\\
{\it University of Wroc\l aw}\\
{\it Plac Maxa Borna 9, 50-204 Wroc\l aw, Poland}\\
\endauthor
\leftheadtext{P. Kosi\'nski, P. Ma\'slanka, J. Sobczyk}
\thanks
*\ \ \ Supported by KBN grant  2 P 302\,217\,06\,p\,02
\endthanks
\thanks
**\ \ \ Supported by KBN grant  2 P 302\,867\,06
\endthanks

\abstract The \2 \4 \3 on the four-\5 $\kappa$-Poincar\'e group and the
corresponding Lie-\0 like structure are described. The \4 \3 on the $n$-\5
\mink{} space covariant under the action of the \1 is constructed.
\endabstract
\endtopmatter

\document

\head I. Introduction
\endhead

In this note we briefly sketch the construction of the \4 \3 on the \1. We
obtain the corresponding Lie \0 structure and prove its equivalence to the
\kapo{} \0. We sketch the  construction of the \4 \3 on the $n$-\5 \mink{}
space covariant under the action of the  \1. Full proofs and
discussions of the further \9 of these \4 calculi will be published
elsewhere.

\head II. The bicovariant calculus on the $\kappa$-Poincar\'e group
\endhead

The \1 $\kap$ is the Hopf *-algebra defined as follows [1]. Consider the
universal  *-algebra with unity generated by the self-adjoint elements
$\vl^\mu{}_\nu$, $x^\mu$ subject to the  following relations:
$$
\aligned
&[x^\mu,x^\nu]  = \dfrac{i}{\kappa}( \de_0{}^\mu x^\nu -  \de_0{}^\nu
x^\mu),\\
&[ \vl^\mu{}_\nu, x^\rho]  = - \dfrac{i}{\kappa} ((\vl^\mu{}_0 -
\de_0{}^\mu) \vl^\rho{}_\nu + (\vl^0{}_\nu - \de_\nu{}^0)  g^{\mu\rho})
\endaligned
\tag{1}
$$
here $g_{\mu\rho} = g^{\mu\rho} = $ diag\,$(+, -, -, -)$ is the  metric
tensor.

The  comultiplication, antipode and counit are introduced as follows:
$$
\aligned
&\vd(\vl^\mu{}_\nu)  =  \vl^\mu{}_\rho \otimes \vl^\rho\nu,\\
&\vd(x^\mu)  = \vl^\mu{}_\nu \otimes v^\nu + x^\mu \otimes I,\\
&S(\vl^\mu{}_\nu)  =  \vl_\nu{}^\mu; \qquad S(x^\mu) = - \vl_\nu{}^\mu
x^\nu,\\
&\ve(\vl^\mu{}_\nu)  =  \de^\mu{}_\nu ; \qquad  \ve(x^\mu) = 0.
\endaligned
\tag{2}
$$
The starting point in our construction of the  \2 *-calculi on the four \5 \1
is
the Wo\-ro\-no\-wicz theory of \4 calculi on  \6 groups [2]. The main
ingredient of this approach is the choice of a right ideal in $\ker \ve$,
which is invariant under the adjoint action of the group. The adjoint action
is defined as follows:
$$
ad(a) = \sum_k b_k \otimes S(a_k)c_k
\tag{3}
$$
here
$$
(\vd \otimes  I) \circ \vd(a) = (I \otimes \vd) \circ \vd(a) = \sum_k a_k
\otimes b_k \otimes c_k.
$$
In the classical case,  the ideal under consideration is $(\ker \ve)^2$. In
order to obtain as slight  a deformation of the classical \3 as possible, we
start
with the generators of $(\ker \ve)^2$. However, it appears that they do not
form a multiplet under the adjoint action of the  \1. To cure this, we modify
them by adding the appropriate $\kappa$-dependent terms. Due to
noncommutativity, the ideal generated in this way is identical with the
whole
$\ker \ve$. Therefore it appears necessary to subtract from the set
spanned by the new generators some ad-invariant  terms (as a consequence the
resulting \3 contains more invariant forms than its classical counterpart).
Finally, we arrive at the following

\proclaim{Theorem 1} Let  $\car \subset \ker \ve$ be the right ideal
generated by the following elements:
$$
\aligned
& (\vl^\a{}_\be - \de^\a{}_\be)(\vl^\mu{}_\nu - \de^\mu{}_\nu); \ \
\widetilde\vd^{\mu\nu\a} \equiv \vd^{\mu\nu\a} - \dfrac{1}{6}
\ve^{\mu\nu\a\be} \ve_{\be\rho\si\de} \vd^{\rho\si\de},\\
& \xti^{\mu\nu} =  x^{\mu\nu} - \dfrac{1}{4} g^{\mu\nu} x^\a{}_\a
\endaligned
\tag{4}
$$
where
$$
\align
&\vd^{\mu\nu\a} = x^\a (\vl^\mu{}_\nu - \de^\mu{}_\nu) - \jak[\de^0{}_\nu
(\vl^{\mu\a} -  g^{\mu\nu}) + \de^\mu{}_0  (\vl^\a{}_\nu - \de^\a{}_\nu)],\\
&  x^{\mu\nu} \equiv  x^{\mu} x^{\nu} + \jak( g^{\mu\a}x^0 -  g^{0\mu}
x^\nu).
\endalign
$$
Then $\car$ has the following \9:
\roster
\item"{(i)}" $\car$ is ad-invariant, ad$(\car) \subset \car \otimes \kap$,
\item"{(ii)}" for any $a \in \car$, $S(a)^* \in \car$,
\item"{(iii)}" $\ker \ve\slash \car$ is spanned by the following elements:
$$
\aligned
& x^\mu; \qquad \vl^\mu{}_\nu - \de^\mu{}_\nu, \ \ \mu < \nu; \qquad \vf
\equiv x^\a{}_\a = x^2 + \dfrac{3i}{\kappa} x^0,\\
& \vf_\mu \equiv \ve_{\mu\nu\a\be} \vd^{\nu\a\be}.
\endaligned
\tag{5}
$$
\endroster
\endproclaim

It is easy to conclude from (iii) that our \3 is fifteen-\5. Having
established the  structure of $\car$, we can now follow closely the
Wo\-ro\-nowicz construction. The basis of the space of the \linv{} 1-forms
consists of the  following elements:
$$
\aligned
& \om^\mu{}_\nu \equiv \pi r^{-1} [I \otimes (\vl^\mu{}_\nu -
\de^\mu{}_\nu)] = \vl_\a{}^\mu d\vl^\a{}_\nu,\\
& \om^\mu \equiv \pi r^{-1} [I \otimes x^\mu] =  \vl_\a{}^\mu
dx^\a,\\
& \om \equiv \pi r^{-1} [I \otimes \vf] = d\vf - 2x_\mu dx^\mu,\\
& \vo_\mu \equiv \pi r^{-1} [I \otimes \vf_\mu] = \ve_{\mu\nu\a\be}
\vl_\si{}^\nu \om^\be \vl^{\si\a} - \dfrac{2i}{\kappa} \ve_{0\mu\a\be}
\om^{\a\be}.
\endaligned
\tag{6}
$$
The next step is to find the commutation rules between the invariant forms
and generators of $\kap$. The detailed calculations result in the following
formulae:
$$
\aligned
& [\vl^\a{}_\be,\om^\mu{}_\nu   ] = 0,\\
& [x^\a,\om^\mu{}_\nu] = -\jak (\de^0{}_\nu \vl^\a{}_\rho \om^{\mu\rho} +
\de^\mu{}_0 \vl^\a{}_\rho \om^\rho{}_\nu - \vl^\a{}_\nu \om^{\mu}{}_0 -
\vl^{\a\mu} \om^0{}_\nu)\\
&\phantom{ [x^\a,\om^\mu{}_\nu] = \ }  - \dfrac{1}{6} {{\ve^\mu}_
\nu}^{\rho\gm}
\vl^\a{}_\rho \vo_\gm,\\
&       [\vl^\mu{}_\nu, \om^\a] = - \jak(\de^0{}_\nu \vl^\mu{}_\rho
\om^{\rho\a} +  \vl^\mu{}_0 \om^\a{}_\nu) -  \dfrac{1}{6} {{\ve^\rho}
_\nu}^{\a\gm} \vl^\mu{}_\rho \vo_\gm,\\
& [x^\mu, \om^\a]  = - \dfrac{1}{4} \vl^{\mu\a} \om + \jak(\vl^{\mu\a}
\om^0 - \de_0{}^\a \vl^\mu{}_\be \om^\be),\\
&       [\vl^\mu{}_\nu,\om]  = \dfrac{4}{\kappa^2} \vl^\mu{}_\rho
\om^\rho{}_\nu,\\
& [x^\mu, \om]  =  \dfrac{4}{\kappa^2} \vl^\mu{}_\rho \om^\rho,\\
& [\vl^\a{}_\be,\vo_\mu] = 0,\\
&[x^\a,\vo_\mu]  = \dfrac{3}{\kappa^2} \ve_{\mu\be\rho\tau}   \vl^{\a\be}
\om^{\rho\tau} - \jak \de^0{}_\mu  \vl^{\a\be} \vo_\be + \jak \vl^{\a}{}_\mu
\vo_0.
\endaligned
\tag{7}
$$
Then, following Woronowicz's paper [2], we can construct the \rinv{} forms:
$$
\aligned
&\eta^\mu{}_\nu  = \om^\be{}_\gm  \vl^\mu{}_\be   \vl_\nu{}^\gm,\\
&\eta^\mu   = - \om^\be{}_\gm  \vl_\rho{}^\gm  x^\rho   \vl^\mu{}_\be +
\om^\be \vl^\mu{}_\be ,\\
&\eta   = \om  ,\\
&\vt_\mu  = \vo_\nu\vl_\mu{}^\nu
\endaligned
\tag{8}
$$
This concludes the description of the bimodule $\vg$ of 1-forms on $\kap$.
The external \0 can  now be constructed as follows [2]. On $\vg^{\otimes 2}$ we
define a bimodule homomorphism $\si$ such that
$$
\si(\om \otimes_{\kap} \eta) = \eta  \otimes_{\kap} \om
\tag{9}
$$
for any \linv{} $\om \in \vg$ and any \rinv{} $\eta \in \vg$. Then, by
definition,
$$
\vg^{\wedge 2} = \dfrac{\vg^{\otimes 2}}{\ker (I - \si)}.
\tag{10}
$$
Higher external powers of $\vg$ can be constructed in a similar way [2].
Eqs. (8)--(10) allow us to calculate the external product of \linv{}
1-forms. The results read:
$$
\aligned
& \om \wedge \om = 0,\\
& \om^\mu{}_\nu \wedge \om^\a{}_\be +  \om^\a{}_\be  \wedge \om^\mu{}_\nu
= 0,\\
& \vo_\a \wedge  \om^\mu{}_\nu + \om^\mu{}_\nu \wedge  \vo\a = 0,\\
& \vo_\mu \wedge  \vo_\nu +  \vo_\nu \wedge   \vo_\mu = 0,\\
& \om^\mu{}_\nu \wedge \om + \om \wedge \om^\mu{}_\nu - \dfrac{4}{\kappa^2}
\om^\si{}_\nu \wedge \om_\si{}^\mu = 0,\\
& \om \wedge \om^\mu + \om^\mu  \wedge \om - \dfrac{4}{\kappa^2}
\om^\mu{}_\si \wedge \om^\si = 0,\\
& \om^\mu \wedge \om^\nu + \om^\nu  \wedge \om^\mu + \jak(\de_0{}^\nu
\om^\mu{}_\rho  \wedge   \om^\rho + \de_0{}^\mu  \om^\nu{}_\rho \wedge
\om^\rho ) = 0,\\
&   \om^\a \wedge \om^{\mu\nu} + \om^{\mu\nu} \wedge\om^\a + \jak(\de_0{}^\nu
\om^{\mu\si} \wedge   \om_\si{}^\a + \de_0{}^\mu  \om^{\si\nu} \wedge
\om_\si{}^\a \\
& \qquad \qquad + \om^{\mu 0} \wedge   \om^{\a\nu}  + \om^{0\nu} \wedge
\om^{\a\mu}) -   \dfrac{1}{6} \ve^{\a\mu\nu\be} X_\be = 0,\\
& \vo_\mu \wedge \om + \om \wedge \vo_\mu - \dfrac{4}{\kappa^2}  \vo_\si
\wedge \om^\si{}_\mu - \dfrac{4}{\kappa^2} X_\mu = 0,\\
& \om^\mu \wedge \vo^\a + \vo^\a  \wedge \om^\mu + \jak (\de_0{}^\a \vo_\rho
\wedge  \om^{\rho\mu} + \vo_0  \wedge \om^{\mu\a})\\
& \qquad \quad + \dfrac{1}{12}  \ve^{\a\mu\rho\si} \vo_\rho \wedge  \vo_\si
- \dfrac{3}{2\kappa^2}    \ve^{\a\si\lm\tau} \om^\mu{}_\si \wedge
\om_{\lm\tau} - \dfrac{1}{4} g^{\mu\a} Y = 0.
\endaligned
\tag{11}
$$
The following notation has been introduced in the above formulae:
$$
\aligned
X_\mu & \equiv \ve_{\mu\rho\si\de} [\om^\rho \wedge \om^{\si\de} +
\om^{\si\de} \wedge \om^\rho + \jak (\de_0{}^\de \om^{\si\lm} \wedge
\om_\lm{}^\rho\\
& \qquad\qquad + \de_0{}^\rho \om^{\lm\de} \wedge \om_\lm + \om^{\si 0}
\wedge \om^{\rho\de} + \om^{ 0\de} \wedge \om^{\rho\si})],\\
Y & \equiv \vo^\rho  \wedge \om_\rho +\om_\rho \wedge \vo^\rho + \jak \vo_\rho
\wedge
   \om^{\rho 0}.
\endaligned
\tag{12}
$$
Eqs. (11) and (12) have serious consequences. The following property of the
classical \3 is usually preserved even in the \6 case: given any basis
$\{\om_i\}$ in the space of (say) \linv{} forms, the basis in $\vg^{\wedge
2}$ is spanned by  $\om_i \wedge \om_j$, $i < j$. This is no longer the case
here. The basis in $\vg^{\wedge 2}$ consists of the following elements:
$$
\align
&\om^{\a\be} \wedge \om^{\mu\nu} \ \ \ ( \a < \be, \ \    \mu < \nu, \ \
(\a \be) \ne (\mu\nu), \ \ \a < \mu); \ \ \ \om^{\a\be} \wedge \om^\mu;\\
&  \om^{\a\be} \wedge \om;  \om^{\a\be} \wedge \vo^{\mu}, \ \ \om^\a \wedge
\om^\mu \ \ (\a < \mu); \ \ \ \om^\a  \wedge \om;\\
& \om^\a \wedge \vo^\mu; \ \ \  \om \wedge \vo^{\mu}; \ \ \vo^\a \wedge
\vo^\mu \ \ (\a < \mu); \ \ \ X_\mu; \ \ Y.
\endalign
$$
Thus  there are five more elements than it is generically expected.

To complete our exterior \3, we derive the Cartan--Maurer equations:
$$
\aligned
& d\om^\mu{}_\nu = \om_\si{}^\mu \wedge \om^\si{}_\nu,\\
& d\om^\mu = \om_\si{}^\mu \wedge \om^\si,\\
& d\om = 0,\\
& d\vo_\mu = - X_\mu - \om_\mu{}^\rho \wedge \vo_\rho.
\endaligned
\tag{13}
$$

In order to obtain the counterpart of the classical Lie \0, we introduce the
\linv{} fields. They are defined by the formula
$$
da = \dfrac{1}{2}(\chi_{{\phantom{}}_{\mu\nu}} \ast a) \om^{\mu\nu} +
(\chi_{{\phantom{}}_{\mu}} \ast
a)\om^\mu + (\chi_{{\phantom{}}_{}} \ast a) \om  + (\lm_{\mu} \ast a)  \vo^\mu
\tag{14}
$$
where, for any linear functional $\vf$ on $\kap$,
$$
\vf \ast a \equiv (I \otimes \vf) \vd(a).
\tag{15}
$$

The product of two functional $\vf_1$, $\vf_2$  is defined by the standard
duality relation
$$
\vf_1 \vf_2(a) \equiv (\vf_1 \otimes \vf_2) \vd(a).
\tag{16}
$$
Finally, we apply the external derivative to both sides of eq. (14). Using
the fact that $d^2a = 0$ and nullifying the coefficients in front of basis
elements of $\vg^{\wedge 2}$,  we find the  \6 Lie \0
$$
\alignat 1
& X^\mu : \ \lm_\mu\Big(1 - \dfrac{4}{\kappa^2} \chi\Big) = - \dfrac{1}{12}
\ve_\mu{}^{\a\rho\si} \chi_{{\phantom{}}_{\a}}
\chi_{{\phantom{}}_{\rho\si}},\tag{17a}\\
& Y  : \ \lm^\mu\chi_{{\phantom{}}_{\mu}} = 0,
\tag{17b}\\
&\om^{\a\be} \wedge \om : \ [\chi_{{\phantom{}}_{\a\be}},\chi] = 0,
\tag{17c}\\
& \om^{\a\be} \wedge \vo^\mu : \ [\chi_{{\phantom{}}_{\a\be}},\lm_\mu]  =
\Big(1 -  \dfrac{4}{\kappa^2} \chi\Big) (g_{\be\mu} \lm_\a - g_{\a\mu}
\lm_\be)\\
& \qquad \qquad \qquad \qquad \qquad \qquad+ \jak \lm_0(g_{\mu\be}
\chi_{{\phantom{}}_{\a}} - g_{\mu\a} \chi_{{\phantom{}}_{\be}})\tag{17d} \\
& \qquad \qquad \qquad \qquad \qquad \qquad + \jak
\de^0{}_\mu(\lm_\a\chi_{{\phantom{}}_{\be}} -
\lm_\be\chi_{{\phantom{}}_{\a}})\\
&\om^{\a} \wedge \om^\mu : \ [\chi_{{\phantom{}}_{\a}},
\chi_{{\phantom{}}_{\mu}}] = 0, \tag{17e}\\
& \om^{\a} \wedge \om : \ [\chi_{{\phantom{}}_{\a}}, \chi] = 0,
\tag{17f}\\
& \om^{\a} \wedge \vo^\mu : \ [\chi_{{\phantom{}}_{\a}}, \lm_{\mu}] = 0,
\tag{17g}\\
& \om  \wedge \vo^\mu : \ [\lm, \lm_\mu] = 0,
\tag{17h}\\
& \vo^\a \wedge \vo^\mu : \ [\lm_\a, \lm_\mu] = \dfrac{1}{6}
\ve_{\a\mu}{}^{\rho\si} \lm_\rho \chi_{{\phantom{}}_{\si}}.
\tag{17i}
\endalignat
$$
In order to simplify the remaining commutation relations we introduce the
following notation:
$$
\aligned
l_i & = \chi_{{\phantom{}}_{i0}},\\
m_i & = \dfrac{1}{2} \ve_{ijk}\chi_{{\phantom{}}_{jk}}.
\endaligned
\tag{18}
$$
They read:\newline
\noindent $\om^{\a\be} \wedge \om^{\mu}$:
$$
[m_i,\chi_{{\phantom{}}_{0}}] = 0,
\tag{19a}
$$
$$
[m_i,\chi_{{\phantom{}}_{k}}] = \Big(1 + \jak \chi_{{\phantom{}}_{0}} -
\dfrac{4}{\kappa^2} \chi\Big)
\ve_{ikl}\chi_{{\phantom{}}_{l}},
\tag{19b}
$$
$$
[l_i,\chi_{{\phantom{}}_{0}}] = \Big(1 + \jak \chi_{{\phantom{}}_{0}} -
\dfrac{4}{\kappa^2}
\chi\Big)\chi_{{\phantom{}}_{i}},
\tag{19c}
$$
$$
[l_i,\chi_{{\phantom{}}_{k}}] = \Big(1 + \jak \chi_{{\phantom{}}_{0}}  -
\dfrac{4}{\kappa^2} \chi\Big)
\de_{ik} \chi_{{\phantom{}}_{0}} ,
\tag{19d}
$$
\noindent $\om^{\a\be} \wedge \om^{\mu\nu}$:
$$
[m_i,m_j] = \Big(1 - \dfrac{4}{\kappa^2} \chi\Big) \ve_{ijk}m_k +
\jak(\chi_{{\phantom{}}_{j}} l_i - \chi_{{\phantom{}}_{i}} l_j) -
\dfrac{6}{\kappa^2}\lm_0
\ve_{ijk}\chi_{{\phantom{}}_{k}},
\tag{19e}
$$
$$
\aligned
[m_i,l_k] & =  \Big(1 + \jak \chi_{{\phantom{}}_{0}} - \dfrac{4}{\kappa^2}
\chi\Big)
\ve_{ikj} l_j + \jak \de_{ik} \chi_{{\phantom{}}_{j}}m_j \\
& - \jak \chi_{{\phantom{}}_{k}}m_i - \dfrac{3}{\kappa^2} (\lm_i
\chi_{{\phantom{}}_{k}} - \de_{ik} \lm_0
\chi_{{\phantom{}}_{0}} - \de_{ik}   \lm_j \chi_{{\phantom{}}_{j}}),
\endaligned
\tag{19f}
$$
$$
[l_i,l_k] = - \Big(1 + \dfrac{2i}{\kappa} \chi_{{\phantom{}}_{0}} -
\dfrac{4}{\kappa^2} \chi
\Big)  \ve_{ikj} m_j -  \dfrac{6}{\kappa^2} \chi_{{\phantom{}}_{0}}
\ve_{ikj}\lm_j.
\tag{19g}
$$
In the $\kappa \to \infty$ limit, the classical Poincar\'e \0 for
$\chi_{{\phantom{}}_{\a\be}}$, $\chi_{{\phantom{}}_{\mu}}$ is restored while
$\lm_\mu$ becomes
proportional to the Pauli--Lubanski four--vector. It can be  checked that
$\chi$
is in turn proportional to the mass squared Casimir operator. Let us note
that the existence of the additional basis elements $X^\mu$, $Y$ of
$\vg^{\wedge 2}$ results in additional relations (17a) and (17b) expressing
(in the  $\kappa \to \infty$  limit) the  Pauli--Lubanski four--vector in
terms
of other generators and the orthogonality of   Pauli--Lubanski and momentum
four--vectors. It is interesting to observe that the analogous relation for
the
mass squared operator in terms of momenta is not derivable in this way
although its validity can be checked. Relations (17d), (17g) and  (17i)
express (in the  $\kappa \to \infty$ limit) the commutation rules for the
Pauli--Lubanski four--vector.

Having our \3 constructed, we can now pose the question: what is the relation
(if any) between our functionals $\chi_{{\phantom{}}_{\mu\nu}}$,
$\chi_{{\phantom{}}_{\a}}$ and  elements of the $\kappa$- Poincar\'e \0
$\kapti$ [3]? It has been shown
recently [4] that $\kap$ and $\kapti$ are formally dual. Therefore we
expect  $\chi_{{\phantom{}}_{\a\be}}$, $\chi_{{\phantom{}}_{\mu}}$, $\lm_\mu$
to be expressible in terms
of elements of $\kapti$. Using, on the one hand, the \9 of the  \linv{}
fields described by Wo\-ro\-nowicz and, on the other hand, the duality
relations $\kap \Longleftrightarrow\kapti$ established in [4], one can prove
that the following substitutions reproduce the \0 and co\0 structure of our \6
Lie \0:
$$
\aligned
& \lm_0 = \dfrac{1}{6} P_iM_i e^{\maks}  ,\\
& \lm_i = \dfrac{1}{6} \Big[  \Big( \kappa \sh\mak + \kub  e^{\maks} \Big)
M_i + \ve_{ijk} P_jN_k e^{\maks} \Big],\\
& m_i = - i e^{\maks} M_i +  \dfrac{6i}{\kappa} \lm_i   ,\\
& l_i = - i e^{\maks} N_i   ,\\
& \chi_{{\phantom{}}_{0}} = -i  \Big( \kappa \sh\mak + \kub  e^{\maks} \Big)
,\\
& \chi_{{\phantom{}}_{i}} = - i P_i e^{\maks} ,\\
& \chi = -\dfrac{1}{8} \Big(2 \kappa^2 \Big(\ch\mak - 1\Big) - {\vec{P\,}^2}
e^{\maks} \Big),
\endaligned
\tag{20}
$$
where $P_\mu$, $M_i$, $N_i$ are generators of the \kapo{} \0.

\head III. The covariant \4 calculi{} on Minkowski spaces
\endhead

The $n$-\5 \1 is defined, [5], by relations (1), (2) where $g_{\mu\nu} =$
diag\,$(+,-,\ld,-)$ is an $n \times n$ matrix. Let us introduce the  $n$-\5
\mink{} space $\cam_\kappa$ as the universal $\ast$-\0 with unity, generated by
selfadjoint elements $y^\mu$ subject to the following relations:
$$
[y^\mu, y^\nu] = \jak(\de_0{}^\mu y^\nu   - \de_0{}^\nu y^\mu).
\tag{21}
$$
 $\cam_\kappa$ can be equipped with the structure of the  \6 group by
defining:
$$
\aligned
\vd  y^\mu = y^\mu \otimes I + I \otimes  y^\mu,\qquad
S(y^\mu) = - y^\mu; \qquad \ve(y^\mu) = 0.
\endaligned
\tag{22}
$$
Let us define the \2 (with respect to the  group structure on
$\cam_\kappa$) \4 \3 on  $\cam_\kappa$. To this end we choose $\car \subset
\ker \ve$ to be the right ideal generated by
$$
\yti^{\mu\nu} \equiv y^\mu y^\nu + \jak(g^{\mu\nu} y^0 - g^{0\mu} y^\nu ) -
\dfrac{1}{n} g^{\mu\nu} \Big( y^2 + \dfrac{(n-1)}{\kappa} iy^0\Big).
\tag{23}
$$
Then we have

\proclaim{Theorem 2} 1)\ $\car$ is $ad$-invariant,  $ad(\car) \subset
\cam_\kappa \otimes \car$; \ 2)\ if $a \in \car$, then $S(a)^* \in  \car$; \
3)\  $\ker \ve\slash \car$ is spanned by $y^\mu$ and $\vf \equiv y^2 +
\frac{n-1}{\kappa} iy^0 $.
\endproclaim

The \linv{} forms are:
$$
\aligned
\tau^\mu & =  \pi r^{-1} (I  \otimes y^\mu) = d y^\mu,\\
\tau & = \pi r^{-1} (I  \otimes \vf) = d\vf - 2 y_\mu d y^\mu.
\endaligned
\tag{24}
$$
They appear to be  \rinv{}, as well. This fact and the  commutation
rules
$$
\aligned
[\tau^\mu, y^\nu]  & = \jak g^{0\mu} \tau^\nu  - \jak g^{\mu\nu} \tau^0 +
\dfrac{1}{n}  g^{\mu\nu
} \tau,\\
[\tau, y^\mu]  & = - \dfrac{n}{\kappa^2}  \tau^\mu
\endaligned
\tag{25}
$$
imply at once the structure of the exterior \3
$$
\aligned
\tau^\mu \wedge \tau^\nu +  \tau^\nu  \wedge \tau^\mu = 0, \qquad
\tau^\mu \wedge \tau + \tau \wedge \tau^\mu = 0.
\endaligned
\tag{26}
$$
Let us proceed to the problem whether the \3 obtained is covariant under
the  action $\rho$ of $\kap$ on $\cam_\kappa$:
$$
\aligned
&\rho(I) = I \otimes I,\\
&\rho(y^\mu) = \vl^\mu{}_\nu \otimes y^\nu + x^\mu \otimes I
\endaligned
\tag{27}
$$
extended by linearity and multiplicativity. Obviously, $\rho$ is a
covariant  action  on  $\cam_\kappa$, i.e. a homomorphism $\rho :
\cam_\kappa \to \kap \otimes \cam_\kappa$ satisfying:
$$
\aligned
&(I \otimes \rho) \circ \rho = (\vd  \otimes I) \circ \rho ,\\
&(\ve  \otimes I)   \circ \rho  = I.
\endaligned
\tag{28}
$$
In order to define the covariant action of $\kap$ on the space of \4
forms on $\cam_\kappa$, we  extend $\rho$ to the action $\rhoti :
\cam_\kappa  \otimes \cam_\kappa \to \kap \otimes \cam_\kappa  \otimes
\cam_\kappa$:
$$
\rhoti (q) = \sum_{i,j,k} a_i{}^kb_i{}^j \otimes x_i{}^k \otimes y_i{}^j
\tag{29}
$$
where $q = \sum_i x_i \otimes y_i \in \cam_\kappa  \otimes \cam_\kappa $ and
$\rho(x_i) = \sum_k a_i{}^k \otimes x_i{}^k$, $\rho(y_i) = \sum_j b_i{}^j
\otimes y_i{}^j$.

Assume now that $\can \in \cam_\kappa^2$ is a sub-bimodule such that
$\rhoti(\can) \subset \kap  \otimes \can$. Then the \4 \3 $(\vg, d)$ defined
by $\can$ has the following property:
$$
\sum_k  x_k dy_k = 0 \rightarrow \sum_k \rho(x_k) (I \otimes d) \rho(y_k) =
0.
\tag{30}
$$
Therefore $\rhoti_1 : \vg \to \kap \otimes \vg$ given by
$$
\rhoti_1 \big(\sum_k  x_k dy_k \big) = \sum_k \rho(x_k) (I \otimes d)
\rho (y_k)
\tag{31}
$$
is a well-defined linear mapping from $\vg $ into $\kap \otimes \vg$.

It  then follows that,  in order to check whether a calculus on $\cam_\kappa$
is
consistent with the action of $\kap $ on  $\cam_\kappa$,   it is sufficient to
check the property $ \rhoti(\can) \subset  \kap \otimes \can$.

One can prove that our first  order \3 is covariant. In a  similar way we
can prove that the higher order \3 is also covariant. We conclude that our
\2 \3 on  $\cam_\kappa$, defined by eqs. (24)--(26), is covariant under the
action of $\kap$ which reads (cf. eq. (31)) as follows:
$$
\aligned
&\rhoti_1 (\tau^\mu) = \vl^\mu{}_\nu \otimes \tau^\nu,\\
&\rhoti_1 (\tau) = I \otimes \tau.
\endaligned
\tag{32}
$$

\enddocument